\documentclass[a4paper, 12pt]{article}
\usepackage[top=1in,bottom=1in,left=1in,right=1in]{geometry}
\usepackage{amsmath}
\usepackage{amsthm}
\usepackage{amsfonts}

\numberwithin{equation}{section}
\newtheorem{thm}{Theorem}[section]
\theoremstyle{definition}
\newtheorem{dfn}{Definition}[section]
\theoremstyle{plain}

\newtheorem{prp}{Proposition}[section]

\begin{document}

\title{\textbf{Kropina change of a Finsler space with m-th root metric}}
\author{\textbf{Gauree Shanker* and Vijeta Singh**}\\\\ *Department of Mathematics and Statistics\\Centaral University of Punjab, Bathinda\\Punjab-151 001,INDIA\\Email-gshankar@cup.ac.in\\\\**Department of Mathematics and Statistics\\Banasthali University,Banasthali\\Rajasthan-304022,INDIA\\Email-90.vijeta@gmail.com}

\date{}
\maketitle

\begin{abstract}
In this paper, we find a condition under which a Finsler space with Kropina change of mth-root metric is projectively related to a mth-root metric and also we find a condition under which this Kropina transformed mth-root metric is locally dually flat. Moreover we find the condition for its Projective flatness.
\end{abstract}
$\newline
$\textbf{Mathematics Subject Classification}: 53B40, 53C60.\newline
$\mathbf{Keywords}$: Finsler space, Kropina space, mth-root metric, projectively related metric, locally dually flat metric, projective flatness.
\section{Introduction.}
Finsler metrics were introduced in order to generalize the Riemannian ones in the sense that the metric should depend not only on the point, but also on the direction. In Finsler geometry, $(\alpha,\beta)$-metrics form a very important and rich class of Finsler metrics which can be expressed in the form $F=\alpha\phi(s),s=\dfrac{\beta}{\alpha}$, where $\alpha$ is a Riemannian metric, $\beta$ is a 1-form and $\phi$ is a positive $C^{\infty}$ function on the domain of definition. In particular, when $\phi=\dfrac{1}{s}$, the Finsler metric $F=\dfrac{\alpha^2}{\beta}$ is called Kropina metric. Kropina metric was first introduced by L. Berwald in connection with a two-dimensional Finsler space with rectilinear extremal and were investigated by V. K. Kropina \cite{7}. They, together with Randers metrics are C-reducible \cite{8}. However, Randers metrics are regular Finsler metrics but Kropina metrics are non-regular Finsler metrics. Kropina metrics seem to be among the simplest nontrivial Finsler metrics with many interesting application in physics, electron optics with a magnetic field, dissipative mechanics and irreversible thermodynamics [\cite{6}, \cite{1}].\\
The theory of mth-root metric, introduced by H.Shimada \cite{10}, and introduced to ecology by Antonelli \cite{3}, had been studied and applied by several authors (\cite{12}, \cite{11} and \cite{19}). It is considered as a generalization of Riemmanian metric in the sense that the second root metric is a Riemmanian metric. For m=3, it is called cubic Finsler metric \cite{8} and for m=4, it is quatric metric \cite{18}. In four-dimension, the special fourth root metric in the form $F=\sqrt[4]{y^1y^2y^3y^4}$ is called Berwald-Moore metric \cite{4}, which is considered as an important subject for a possible model of space time by physicists. In this paper,  we find the condition under which the transformed Finsler space is projectively related with given Finsler space. Also, we find the condition under which the transformed Finsler space is locally dually flat. Moreover, we find the condition for its Projectively flatness.
\section{Preliminaries.}
 Let $M^n$ be a real smooth manifold of dimension $n$, $T_xM$ denotes the tangent space of $M^n$ at $x$. The tangent bundle $TM$ is the union of tangent spaces, $TM:=\bigcup_{x\in M}T_xM $.We denote the elements of $TM$ by $(x,y)$, where $x=(x^i)$ is a point of $M^n$ and $y\in T_xM$ called supporting element. We denote $TM_0=TM\setminus\{0\}$.
\begin{dfn}
A Finsler metric on a manifold $M$ is a $C^{\infty}$ function $F: TM\backslash\{0\}\rightarrow[0,\infty)$ satisfying the following conditions:\\
(1) Regularity: $F$ is $C^{\infty}$ on $TM\backslash\{0\}$ .\\
(2) Positive homogeneity: $F(x,\lambda y) = \lambda F(x,y), \lambda>0$.\\
(3) Strong convexity: the fundamental tensor $g_{ij}(x,y)$ is positive definite for all $(x,y)\in TM\backslash\{0\}$ , where $g_{ij}=\dfrac{1}{2}\dfrac{\partial^2F^2}{\partial y^i\partial y^j}$.\\
\end{dfn}
By the homogeneity of $F$, we have $F(x,y)=\sqrt{g_{ij}(x,y)y^iy^j}$. An important class of Finsler metrics are Riemann metrics, which are in the form of $F(x,y)=\sqrt{g_{ij}(x)y^iy^j}$. Another important class of Finsler metrics are Minkowski metrics, which are in the form of $F(x,y)=\sqrt{g_{ij}(y)y^iy^j}$.\\
The pair $(M^n,F)=F^n$ is called a Finsler space, F is called the fundamental function and $g_{ij}$ is called the fundamental tensor of the Finsler space $F^n$. \\
The normalized supporting element $l_i$, angular metric tensor $h_{ij}$ and metric tensor $g_{ij}$ of $F^n$ are defined respectively as:
\begin{align}
l_i=\dfrac{\partial F}{\partial y^i}, && h_{ij}=F\dfrac{\partial^2F}{\partial y^i\partial y^j}, && g_{ij}=\dfrac{1}{2}\dfrac{\partial^2F^2}{\partial y^i\partial y^j}.
\end{align}
Let $F$ be a Finsler metric defined by $F=\sqrt[m]{A}$, where $A$ is given by $A:=a_{i_1i_2...i_m}(x)y^{i_1}y^{i_2}...y^{i_m}$, with $a_{i_1...i_m}$ symmetric in all its indices \cite{10}. Then $F$ is called an mth-root Finsler metric. Clearly, $A$ is homogeneous of degree $m$ in $y$.\\
Let 
\begin{equation}
A_i=a_{ii_2...i_m}(x)y^{i_2}...y^{i_m}=\dfrac{1}{m}\dfrac{\partial A}{\partial y^i},
\end{equation}
\begin{equation}
A_{ij}=a_{iji_3...i_m}(x)y^{i_3}...y^{i_m}=\dfrac{1}{m(m-1)}\dfrac{\partial^2A}{\partial y^i\partial y^j},
\end{equation}
\begin{equation}
A_{ijk}=a_{ijki_4...i_m}(x)y^{i_4}...y^{i_m}=\dfrac{1}{m(m-1)(m-2)}\dfrac{\partial^3A}{\partial y^i\partial y^j\partial y^k}.
\end{equation}
The normalized supporting element of $F^n$ is given by
\begin{equation}
l_i:=F_{y^i}=\dfrac{\partial F}{\partial y^i}=\dfrac{\partial A^{\dfrac{1}{m}}}{\partial y^i}=\dfrac{1}{m}\dfrac{\dfrac{\partial A}{\partial y^i}}{A^{\dfrac{(m-1)}{m}}}=\dfrac{A_i}{F^{m-1}}.
\end{equation}
Let us consider the transformation
\begin{equation}
\bar{F}=\dfrac{F^2}{\beta},
\end{equation}
where $F=\sqrt[m]{A}$ is an mth-root metric and $\beta=b_i(x)y^i$ is a one form on the manifold $M^n$. Clearly, $\bar{F}$ is also a Finsler metric on $M^n$, given by Kropina change mth root metric. Throughout the paper, we call the Finsler metric $\bar{F}$ as transformed mth root metric and $(M^n,\bar{F})=\bar{F}^n$ as transformed Finsler space. We restrict ourselves for $m>2$ throughout the paper and also the quantities corresponding to the transformed Finsler space $\bar{F}^n$ will be denoted by putting bar on the top of that quantity. 

 \section{Fundamental metric tensor of Kropina transformed Finsler space with mth-root metric}
 Let us consider the Finsler metric given in (2.6), where $F=\sqrt[m]{A}$ and $\beta=b_i(x)y^i$ is a differential one form on the manifold $M^n$. This metric is called Kropina change of Finsler metric.\\
 The differentiation of (2.6) with respect to $y^i$ yields the normalized supporting element $\bar{l_i}$ given by
 \begin{equation}
 \bar{l_i}=\dfrac{2F}{\beta}l_i-\dfrac{F^2}{\beta^2}b_i.
 \end{equation}
 In view of (2.5), we have 
 \begin{equation}
 \bar{l_i}=\dfrac{2A_i}{\beta F^{m-2}}-\dfrac{F^2}{\beta^2}b_i. 
 \end{equation}
 Again differentiation of (3.2) with respect to $y^j$ yields:
 \begin{equation}
 \bar{h_{ij}}=\dfrac{2\bar{F}}{\beta}\bigg[\dfrac{(m-1)}{F^{m-2}}A_{ij}-\dfrac{(A_ib_j+A_jb_i)}{F^{m-2}\beta}+\dfrac{b_ib_j}{F^2\beta^2}-\dfrac{(m-2)}{F^{2(m-1)}}A_iA_j\bigg].
 \end{equation}
 From (3.2) and (3.3), the fundamental metric tensor $\bar{g_{ij}}$ of Finsler space $\bar{F}^n$ is given by:
 \begin{equation*}
 \bar{g_{ij}}=\bar{h_{ij}}+\bar{l_i}\bar{l_j}
 \end{equation*}
 After simplification, we get
 \begin{align}
 \bar{g_{ij}}=2\tau^2\bigg(\dfrac{(m-1)}{F^{(m-2)}}A_{ij}-\dfrac{2\tau}{F^{(m-1)}}(A_ib_j+A_jb_i)+\tau^2(\dfrac{1}{F^4}+\dfrac{1}{2})b_ib_j-\dfrac{(m-4)}{F^{2(m-1)}}A_iA_j\bigg),
 \end{align}
 where $\tau=\dfrac{\bar{F}}{F}=\dfrac{F}{\beta}$.\\
 The contravariant metric tensor $\bar{g^{ij}}$ of Finsler space $\bar{F}^n$ is given by
 \begin{align}
 \bar{g^{ij}}=\dfrac{F^{m-2}}{2\tau^2(m-1)}A^{ij}+p_0b^ib^j+\dfrac{2\beta^3(m-4)p_1}{F^2(m-1)[\beta^4(m-4)-8F^4d^2]}(b^iy^j+b^jy^i)+p_2y^iy^j.
 \end{align}
 where 
 \begin{align}
 \notag&p_0=\dfrac{4F^m[1+q(q(1+v)-(3+v))]}{\beta^2[(m-4)-8\tau^4d^2]};  \delta=\dfrac{-8F^4}{\beta^4(m-4)}\\
 \notag&q=\dfrac{\delta w^2}{1+\delta c^2}; w=\dfrac{F^{m-2}}{2\tau^2(m-1)};c^2=\dfrac{F^{m-3}\beta b^2}{2\tau(m-1)};v=\dfrac{(m-4)\beta}{2F^m};\\
 \notag&d^2=w[v\beta+v^2F^m+(b^2+v\beta)(1-\dfrac{\delta w(1+v)}{1+\delta c^2})];\\
 \notag&p_1=\dfrac{8(m-1)^2\tau^4+2\delta F^{2(m-2)}+\delta F^{m-4}(m-4)\beta}{4\tau^4(m-1)^2+\delta F^{m-2}b^2\tau^4};\\
  &p_2=\dfrac{(m-4)^2}{2F^6[(m-4)-8\tau^4d^2]}.
 \end{align}
 Here we have used $A^{ij}A_j=A^i=y^i$ and $A_jb^j=\beta$.
 
 \begin{prp}
 The Fundamental metric $\bar{g_{ij}}$ and its inverse tensor $\bar{g^{ij}}$ of Kropina transformed mth-root Finsler space $\bar{F}^n$ are given by equation (3.4) and (3.5) respectively.
 \end{prp}
 
 \section{Spray coefficients of the Finsler space given by Kropina change of mth-root metric.}
   The  geodesics of a Finsler space $F^n$ are given by the following system of equations
   \begin{equation*}
   \dfrac{d^2x^i}{dt^2}+G^i\bigg(x,\dfrac{dx}{dt}\bigg)=0,
   \end{equation*}
   where
   \begin{equation}
  G^i=\dfrac{1}{4}g^{il}\{[F^2]_{x^ky^l}y^k-[F^2]_{x^l}\}.
   \end{equation}
    The local functions $G^i=G^i(x,y)$ define a global vector field $G=y^i\dfrac{\partial}{\partial x^i}-2G^i(x,y)\dfrac{\partial}{\partial y^i}$ on $TM$. $G$ is called spray of $F$ and $G^i$ is called spray coefficient.\\
   Two Finsler metrics $F$ and $\bar{F}$ on a manifold $M^n$ are called projectively related if there is a scalar function $P(x,y)$ defined on $TM_0$ such that $\bar{G^i}=G^i+Py^i$, where $\bar{G^i}$ and $G^i$ are the geodesic spray coefficients of $\bar{F^n}$ and $F^n$ respectively. In other words two metrics $\bar{F}$ and $F$ are called projectively related if any geodesic of the first is also good geodesic for the second and vice versa.\\
   In view of equation (3.4) the metric tensor $\bar{g_{ij}}$ of $\bar{F^n}$ can be rewritten as :
   \begin{align}
   \bar{g}_{ij}=2\tau^2g_{ij}-\dfrac{4\tau^3}{F^{m-1}}(A_ib_j+A_jb_i)+\dfrac{(2+F^4)}{\beta^4}b_ib_j+\dfrac{4\tau^2}{F^{2(m-1)}}A_iA_j,
   \end{align}
   where
   \begin{align}
   g_{ij}=(m-1)\dfrac{A_{ij}}{F^{m-2}}-(m-2)\dfrac{A_iA_j}{F^{2(m-1)}}.
   \end{align}
   Further, in view of equation (3.5) contravariant metric tensor $\bar{g^{ij}}$ can be rewritten as :
   \begin{align}
   \bar{g}^{ij}=\dfrac{g^{ij}}{2\tau^2}+p_0b^ib^j+\dfrac{2\beta^3(m-4)p_1}{F^2(m-1)[\beta^4(m-4)-8F^4d^2]}(b^iy^j+b^jy^i)+p_3y^iy^j,
   \end{align}
   where $p_0$, $p_1$ and $d^2$ are expressed in equation(3.6) and
   \begin{align}
  \notag &g^{ij}=\dfrac{F^{m-2}}{(m-1)}A^{ij}+\dfrac{(m-2)y^iy^j}{(m-1)F^2},\\
  &p_3=\dfrac{(m-4)^2(m-1)\tau^2-(m-2)[(m-4)-8\tau^4d^2]\beta^4}{2F^2\tau^2\beta^4(m-1)[(m-4)-8\tau^4d^2]}.
    \end{align}
    The spray coefficients of Kropina transformed Finsler space $F^n$ are given by
    \begin{equation*}
   \bar{G^i}=\dfrac{1}{4}\bar{g}^{il}\{[\bar{F^2}]_{x^ky^l}y^k-[\bar{F^2}]_{x^l}\}
    \end{equation*}
    It can also be written as
    \begin{equation}
    \bar{G}^i=\dfrac{1}{4}\bar{g}^{il}\bigg[\bigg(2\dfrac{\partial\bar{g_{jl}}}{\partial x^k}-\dfrac{\partial\bar{g_{jk}}}{\partial x^l}\bigg)y^jy^k\bigg].
    \end{equation}
    From (4.2), (4.4) and (4.6), we get
    \begin{align*}
    \bar{G}^i&=\dfrac{\bar{g}^{il}}{4}\bigg[\bigg(\dfrac{\partial}{\partial x^k}\{2\tau^2g_{jl}-\dfrac{4\tau^3}{F^{m-1}}(A_jb_l+A_lb_j)+\dfrac{(2+F^4)}{\beta^4}b_lb_j+\dfrac{4\tau^2}{F^{2(m-1)}}A_lA_j\}-\\&\dfrac{\partial}{\partial x^l}\{2\tau^2g_{jk}-\dfrac{4\tau^3}{F^{m-1}}(A_jb_k+A_kb_j)+\dfrac{(2+F^4)}{\beta^4}b_kb_j+\dfrac{4\tau^2}{F^{2(m-1)}}A_kA_j\}\bigg)y^jy^k\bigg]
    \end{align*}
    which implies that
    \begin{align*}
    \bar{G}^i&=\dfrac{\bar{g}^{il}}{4}\bigg[\bigg(2\{2\tau^2\dfrac{\partial g_{jl}}{\partial x^k}+g_{jl}\dfrac{\partial}{\partial x^k}(2\tau^2)+\dfrac{\partial X_{jl}}{\partial x^k}\}-\{2\tau^2\dfrac{\partial g_{kl}}{\partial x^l}+g_{jk}\dfrac{\partial}{\partial x^l}(2\tau^2)+\dfrac{\partial X_{jk}}{\partial x^l}\}\bigg)y^jy^k\bigg],
    \end{align*}
    where
    \begin{equation*}
    X_{jl}=-\dfrac{4\tau^3}{F^{m-1}}(A_jb_l+A_lb_j)+\dfrac{(2+F^4)}{\beta^4}b_lb_j+\dfrac{4\tau^2}{F^{2(m-1)}}A_lA_j.
    \end{equation*}
    Now
    \begin{align*}
   \bar{G}^i=\dfrac{1}{4}&\bigg[\dfrac{g^{il}}{2\tau^2}+y^i(p_1b^l+p_3y^l)+b^i(p_0b^l+p_1y^l)\bigg]\times \\ &\bigg[2\tau^2\bigg(2\dfrac{\partial g_{jl}}{\partial x^k}-\dfrac{\partial g_{jk}}{\partial x^l}\bigg)+2\omega_kg_{jl}-\omega_lg_{jk}+2\dfrac{\partial X_{jl}}{\partial x^k}-\dfrac{\partial X_{jk}}{\partial x^l}\bigg]y^jy^k,
    \end{align*}
    where $\omega_k=\dfrac{\partial}{\partial x^k}(2\tau^2)$.\\\\
    Simplifying, we get
   \begin{align*}
    \bar{G}^i&=\dfrac{g^{il}}{4}\bigg(2\dfrac{\partial g_{jl}}{\partial x^k}-\dfrac{\partial g_{jk}}{\partial x^l}\bigg)y^jy^k+\dfrac{F^{m-2}A^{il}}{8\tau^2(m-1)}[2\omega_kg_{jl}-\omega_lg_{jk}+2\dfrac{\partial X_{jl}}{\partial x^k}-\dfrac{\partial X_{jk}}{\partial x^l}]y^jy^k\\&+\dfrac{(m-2)}{4\tau^2F^2(m-1)}y^iy^l\bigg[2\omega_kg_{jl}-\omega_lg_{jk}+2\dfrac{\partial X_{jl}}{\partial x^k}-\dfrac{\partial X_{jk}}{\partial x^l}\bigg]+\dfrac{1}{4}\bigg[y^i(p_1b^l+p_3y^l)+b^i(p_0b^l+p_1y^l)\bigg]\times \\&\bigg[2\tau^2\bigg(\dfrac{\partial g_{jl}}{\partial x^k}-\dfrac{\partial g_{jk}}{\partial x^l}\bigg)+2\omega_kg_{jl}-\omega_lg_{jk}+2\dfrac{\partial X_{jl}}{\partial x^k}-\dfrac{\partial X_{jk}}{\partial x^l}\bigg]y^jy^k.
    \end{align*}
    The above equation may be written as
    \begin{equation*}
    \bar{G}^i=G^i+Py^i+Q^i,
    \end{equation*}
    where
    \begin{align*}
    P&=\dfrac{1}{4}y^i(p_1b^l+p_3y^l)\times\bigg[2\tau^2\bigg(\dfrac{\partial g_{jl}}{\partial x^k}-\dfrac{\partial g_{jk}}{\partial x^l}\bigg)+2\omega_kg_{jl}-\omega_lg_{jk}+2\dfrac{\partial X_{jl}}{\partial x^k}-\dfrac{\partial X_{jk}}{\partial x^l}\bigg]y^jy^k\\&+\dfrac{(m-2)}{4\tau^2F^2(m-1)}y^ly^jy^k\bigg[2\omega_kg_{jl}-\omega_lg_{jk}+2\dfrac{\partial X_{jl}}{\partial x^k}-\dfrac{\partial X_{jk}}{\partial x^l}\bigg]
    \end{align*}
    and
    \begin{align*}
    Q^i&=\dfrac{1}{4}b^i(p_2b^l+p_1y^l)\times\bigg[2\tau^2\bigg(\dfrac{\partial g_{jl}}{\partial x^k}-\dfrac{\partial g_{jk}}{\partial x^l}\bigg)+2\omega_kg_{jl}-\omega_lg_{jk}+2\dfrac{\partial X_{jl}}{\partial x^k}-\dfrac{\partial X_{jk}}{\partial x^l}\bigg]y^jy^k\\&+\dfrac{F^{m-2}A^{il}}{8\tau^2(m-1)}y^jy^k\bigg[2\omega_kg_{jl}-\omega_lg_{jk}+2\dfrac{\partial X_{jl}}{\partial x^k}-\dfrac{\partial X_{jk}}{\partial x^l}\bigg].
    \end{align*}
    The metrics $\bar{F}$ and $F$ are projectively related if $Q^i=0$, which implies
    \begin{align}
    \notag&\dfrac{1}{4}b^i(p_2b^l+p_1y^l)\times\bigg[2\tau^2\bigg(\dfrac{\partial g_{jl}}{\partial x^k}-\dfrac{\partial g_{jk}}{\partial x^l}\bigg)+2\omega_kg_{jl}-\omega_lg_{jk}+2\dfrac{\partial X_{jl}}{\partial x^k}-\dfrac{\partial X_{jk}}{\partial x^l}\bigg]y^jy^k\\&=\dfrac{F^{m-2}A^{il}}{8\tau^2(m-1)}y^jy^k\bigg[2\omega_kg_{jl}-\omega_lg_{jk}+2\dfrac{\partial X_{jl}}{\partial x^k}-\dfrac{\partial X_{jk}}{\partial x^l}\bigg].
    \end{align}
    Thus we have the following:
    \begin{thm}
    The Kropina transformed mth-root metric $\bar{F}$ and mth-root metric $F$, on an open subset $\cup \subset R^n$, are projectively related if equation (4.7) is satisfied.
    \end{thm}
    
    \section{Locally dually flatness of a Finsler Space with Kropina changed mth-root metric}
       The notion of dually flat metrics was first introduced by S.I. Amari and H. Nagaoka \cite{1} when they studied the information geometry on Riemannian spaces. Later on, Z. Shen extended the notion of dually flatness to Finsler metrics \cite{19}.\\
       A Kropina changed Finsler space on manifold $M^n$ is said to be locally dually flat if at any point there is a standard coordinate system $(x^i,y^i)$ in $TM$ such that $[\bar{F}^2]_{x^ky^l}y^k=2[\bar{F}^2]_{x^l}$. The coordinate $x^i$ is called an adapted local coordinate system. It is well known that every locally Minkowskian metric is locally flat. \\
       Consider the Kropina changed mth root metric $\bar{F}=\dfrac{F^2}{\beta}$, where $F$ is an mth-root metric.. Then, we have
       \begin{align}
       [\bar{F}^2]_{x^l}=\bigg[\dfrac{F^2}{\beta}\bigg]^2_{x^l}=2\bigg[\dfrac{F^2}{\beta}\bigg]\bigg[\dfrac{\frac{2}{m}FA^{\frac{1}{m}-1}A_{x^l}}{\beta}-\dfrac{F^2\beta_l}{\beta^2}\bigg].
       \end{align}
       If we put $b_{ij}=\dfrac{\partial b_i}{\partial x^j}$, then $\beta_j=\dfrac{\partial \beta}{\partial x^j}=b_{ij}y^j$.\\
       On simplifying equation (5.1) further, we get
     \begin{align}
     [\bar{F}^2]_{x^l}=\dfrac{4}{m\beta^2}A^{\frac{4-m}{m}}A_{x^l}-\dfrac{2}{\beta^3}A^{\frac{4}{m}}\beta_l
     \end{align}
     From (5.2), we get
     \begin{align}
    [\bar{F}^2]_{x^k}=\dfrac{4}{m\beta^2}A^{\frac{4-m}{m}}A_{x^k}-\dfrac{2}{\beta^3}A^{\frac{4}{m}}\beta_k.
     \end{align}
     and 
     \begin{align}
     \notag[\bar{F}^2]_{x^ky^l}&=\dfrac{4}{m\beta^2}A^{\frac{4-m}{m}}A_{x^ky^l}+\dfrac{4}{m\beta^2}\bigg(\dfrac{4-m}{m}\bigg)A^{\frac{4-2m}{m}}A_{y^l}-\dfrac{8}{m\beta^3}A^{\frac{4-m}{m}}A_{x^k}b_l\\&-\dfrac{8}{m\beta^3}A^{\frac{4-m}{m}}A_{y^l}\beta_l+\dfrac{6}{\beta^4}A^{\frac{4}{m}}\beta_kb_l+\dfrac{2}{\beta^3}A^{\frac{4}{m}}b_lk.
     \end{align}
     For the Finsler metric $\bar{F}$ to be locally dually flat, we must have
     \begin{equation}
     [\bar{F}^2]_{x^ky^l}y^k-2[\bar{F}^2]_{x^l}=0.
     \end{equation}
     Therefore, from (5.2), (5.4) and (5.5), we get
     \begin{align*}
    \notag [\bar{F}^2]_{x^ky^l}y^k-2[\bar{F}^2]_{x^l}=&\dfrac{4}{m\beta^2}A^{\frac{4-m}{m}}A_{x^ky^l}y^k+\dfrac{4}{m\beta^2}\bigg(\dfrac{4-m}{m}\bigg)A^{\frac{4-2m}{m}}A_{y^l}y^k-\dfrac{8}{m\beta^3}A^{\frac{4-m}{m}}A_{x^k}b_ly^k\\\notag -&\dfrac{8}{m\beta^3}A^{\frac{4-m}{m}}A_{y^l}\beta_ly^k+\dfrac{6}{\beta^4}A^{\frac{4}{m}}\beta_kb_ly^k+\dfrac{2}{\beta^3}A^{\frac{4}{m}}b_lky^k\\&-\dfrac{8}{m\beta^2}A^{\frac{4-m}{m}}A_{x^l}+\dfrac{4}{\beta^3}A^{\frac{4}{m}}\beta_l=0.
     \end{align*}
     which implies
     \begin{align*}
     \notag A_{x^l}\bigg[\dfrac{8}{m\beta^2}A^{\frac{4-m}{m}}\bigg]&=\dfrac{4}{m\beta^2}\bigg(\dfrac{4-m}{m}\bigg)A^{\frac{4-2m}{m}}A_{yl}A_{0}+\dfrac{4}{m\beta^2}A^{\frac{4-m}{m}}A_{0l}\\\notag &-\dfrac{8}{m\beta^3}A^{\frac{4-m}{m}}A_0b_l-\dfrac{8}{m\beta^3}A^{\frac{4-m}{m}}A_{y^l}y^k+\dfrac{6}{\beta^4}A^{\frac{4}{m}}\beta_kb_ly^k\\ &+\dfrac{2}{\beta^3}A^{\frac{4}{m}}\beta_l-\dfrac{8}{m\beta^2}A^{\frac{4-m}{m}}A_{x^l}+\dfrac{4}{\beta^3}A^{\frac{4}{m}}b_l,
     \end{align*}
     where
     \begin{align*}
     A_0=A_{x^k}y^k, && A_{0l}=A_{x^ky^l}y^k
     \end{align*}
     Therefore $\bar{F}$ is locally dually flat metric iff
     \begin{align}
     A_{x^l}=\dfrac{1}{2\beta F^2}\bigg(\dfrac{4-m}{m}\bigg)A_0A_{yl}+\dfrac{1}{2}A_{0l}-\dfrac{1}{\beta}A_0b_l-\dfrac{1}{\beta}A_{yl}y^k+\dfrac{3m}{4\beta^2}A\beta_ky^kb_l+\dfrac{m}{2\beta}A\beta_l-A_{x^l}+\dfrac{m}{2\beta}Ab_l.
     \end{align}
     Thus we have 
     \begin{thm}
     Let $\bar{F}$ be a Kropina changed mth-root Finsler metric on a Finsler manifold $M^n$. Then, $\bar{F}$ is locally dually flat metric iff (5.6) holds.
     \end{thm}
     
     \section{Projectively Flatness of a Finsler Space with Kropina changed mth-root metric}
     There is another important notion in Finsler geometry, known as projectively projective flatness of Finsler metrics. A Finsler metric $\bar{F}=\bar{F}(x,y)$ on an open subset $\cup \subset R^n$ is projectively flat iff it satisfies following equation:
     \begin{equation*}
     \bar{F}_{x^ky^l}y^k-F_{x^l}=0.
     \end{equation*} 
     Since we have $\bar{F}=\dfrac{F^2}{\beta}$, where $F=\sqrt[m]{A}$, we get      
     \begin{align}
     [\bar{F}]_{x^l}=\dfrac{2}{m\beta}A^{\frac{2-m}{m}}A_{x^l}-\dfrac{A^{\frac{2}{m}}}{\beta^2}\beta_l.
     \end{align}
     From (6.1)
      we get
     \begin{equation*}
     [\bar{F}]_{x^k}=\dfrac{2}{m\beta}A^{\frac{2-m}{m}}A_{x^k}-\dfrac{A^{\frac{2}{m}}}{\beta^2}\beta_k
     \end{equation*}
     which implies
     \begin{align}
     \notag [\bar{F}]_{x^ky^l}=&\dfrac{2}{m\beta}\bigg(\dfrac{2-m}{m}\bigg)A^{\frac{2(1-m)}{m}}A_{y^l}A_{x^k}+\dfrac{2}{m\beta}A^{\frac{2-m}{m}}A_{x^ky^l}-\dfrac{2}{m\beta^2}A^{\frac{2-m}{m}}A_{x^k}b_l\\&-\dfrac{2}{m\beta^2}A^{\frac{2-m}{m}}A_{y^l}\beta_k+\dfrac{2}{\beta^3}A^{\frac{2}{m}}\beta_kb_l-\dfrac{F^2}{\beta^2}b_{lk}.
     \end{align}
     For the Finsler metric $\bar{F}$ to be Projectively Flat, we must have
     \begin{equation}
     \bar{F}_{x^ky^l}y^k-F_{x^l}=0.
     \end{equation}
     Therefore, from (6.1), (6.2) and (6.3), we obtain
     \begin{align}
     \notag \bar{F}_{x^ky^l}y^k-F_{x^l}=&\dfrac{2}{m\beta}\bigg(\dfrac{2-m}{m}\bigg)A^{\frac{2(1-m)}{m}}A_{y^l}A_{x^k}y^k+\dfrac{2}{m\beta}A^{\frac{2-m}{m}}A_{x^ky^l}y^k-\dfrac{2}{m\beta^2}A^{\frac{2-m}{m}}A_{x^k}b_ly^k\\&-\dfrac{2}{m\beta^2}A^{\frac{2-m}{m}}A_{y^l}\beta_ky^k+\dfrac{2}{\beta^3}A^{\frac{2}{m}}\beta_kb_ly^k-\dfrac{F^2}{\beta^2}b_{lk}y^k-\dfrac{2}{m\beta}A^{\frac{2-m}{m}}A_{x^l}+\dfrac{A^{\frac{2}{m}}}{\beta^2}\beta_l=0
     \end{align}
    which implies
     \begin{align}
     \notag A_{x^l}\bigg[\dfrac{2}{m\beta}A^{\frac{2-m}{m}}\bigg]&=\dfrac{2}{m\beta}\bigg(\dfrac{2-m}{m}\bigg)A^{\frac{2(1-m)}{m}}A_{y^l}A_0+\dfrac{2}{m\beta}A^{\frac{2-m}{m}}A_{0l}-\dfrac{2}{m\beta^2}A^{\frac{2-m}{m}}A_0b_l\\&-\dfrac{2}{m\beta^2}A^{\frac{2-m}{m}}A_{y^l}\beta_ky^k+\dfrac{2}{\beta^3}A^{\frac{2}{m}}\beta_kb_ly^k,
     \end{align}
     where
     \begin{align*}
      A_0=A_{x^k}y^k, && A_{0l}=A_{x^ky^l}y^k.
      \end{align*}
       Therefore $\bar{F}$ is projectively flat metric iff
       \begin{align}
       A_{x^l}=\bigg(\dfrac{2-m}{m}\bigg)A^{-1}A_{y^l}A_0+A_{0l}-\dfrac{1}{\beta}A_0b_l-\dfrac{1}{\beta}A_{y^l}\beta_ky^k+\dfrac{m}{\beta}A\beta_kb_ly^k.
       \end{align}
       Thus we have the following: 
        \begin{thm}
        Let $\bar{F}$ be a Kropina changed mth-root Finsler metric on a Finsler manifold $M^n$. Then, $\bar{F}$ is projectively flat metric iff (6.6) holds.
        \end{thm}
        
   \end{document}